\documentclass[12pt,a4paper]{article}
\usepackage{amsmath}
\usepackage{indentfirst}
\usepackage{geometry}
\usepackage{amsfonts}
\usepackage{amssymb}
\usepackage{amsthm}
\usepackage{mathrsfs}
\usepackage{enumerate}

\newtheorem{Assumption}{Assumption}[section]
\newtheorem{Theorem}{Theorem}[section]
\newtheorem{Remark}[Theorem]{Remark}
\newtheorem{Lemma}[Theorem]{Lemma}

\newtheorem{Proposition}[Theorem]{Proposition}

\numberwithin{equation}{section}
\def\C_lip{\mathcal{C}_{l,lip}}

\newcommand{\cF}{{\mathcal F}}

\newcommand{\R}{\mathbb{R}}

\newcommand{\p}{\mathbb{P}}

\newcommand{\E}[1]{\mathbb{E}\left[#1\right]}

\newcommand{\abs}[1]{\left|#1\right|} 
\newcommand{\ep}{\varepsilon} 
\renewcommand{\geq}{\geqslant} 
\renewcommand{\leq}{\leqslant} 
\renewcommand{\ge}{\geqslant} 
\renewcommand{\le}{\leqslant}

\def\be{\mathbb{E}}

\def\setu{\mathcal{U}[0,T]}
\def\setad{\mathcal{U}_{ad}[0,T]}
\def\setk{\mathbb{K}}
\def\trans{\intercal}

\def\setmum{{\mathcal{M}^-}([0,T])}
\def\setmu{{\mathcal{M}^+}([0,T])}

\def\haty{\widehat{Y}}

\usepackage{hyperref}
\hypersetup{
colorlinks=true, 
breaklinks=true, 
urlcolor=blue, 
linkcolor=red, 
citecolor=green, 
pdftitle={}, 
pdfauthor={}, 
pdfsubject={} 
}

\usepackage{ulem}

\title{Optimal control of SDEs with expected path constraints and related constrained FBSDEs\footnotemark[1]}

\author{Ying Hu\footnotemark[2] \and Shanjian Tang\footnotemark[3] \and Zuo Quan Xu\footnotemark[4]}

\begin{document}
\maketitle
\begin{abstract}
In this paper, we consider optimal control of stochastic differential equations subject to an expected path constraint. The stochastic maximum principle is given for a general optimal stochastic control in terms of constrained FBSDEs. In particular, the compensated process in our adjoint equation is deterministic, which seems to be new in the literature. 
For the typical case of linear stochastic systems and quadratic cost functionals (i.e., the so-called LQ optimal stochastic control), a verification theorem is established, and the existence and uniqueness of the constrained reflected FBSDEs are also given.
\end{abstract}

\footnotetext[1]{Ying Hu is partially supported by Lebesgue Center of Mathematics ``Investissements d'avenir" program-ANR-11-LABX-0020-01, by CAESARS-ANR-15-CE05-0024 and by MFG-ANR-16-CE40-0015-01. Shanjian Tang is partially supported by National Science Foundation of China (Nos. 11631004 and 12031009). Zuo Quan Xu is partially supported by NSFC (No.11971409),  The Research Grants Council of Hong Kong (GRF No.15202421), The PolyU-SDU Joint Research Center on Financial Mathematics, The CAS AMSS-POLYU Joint Laboratory of Applied Mathematics, and The Hong Kong Polytechnic University. }
\footnotetext[2]{Univ. Rennes, CNRS, IRMAR-UMR6625, F-35000, Rennes, France; School of Mathematical Sciences, Fudan University, Shanghai 200433,
China (\textit{E-mail: ying.hu@univ-rennes1.fr}).}
\footnotetext[3]{Department of Finance and Control Sciences, School of
Mathematical Sciences, Fudan University, Shanghai 200433, China (\textit{e-mail: sjtang@fudan.edu.cn}).}
\footnotetext[4]{Department of Applied Mathematics, The Hong Kong Polytechnic University, Kowloon, Hong Kong, China (\textit{e-mail: maxu@polyu.edu.hk}).}

{\bf MSC2020 subject classifications:} Primary 93E20, 60H30; secondary 60H10, 93E03. 


\section{Introduction} 

In this paper, we consider the following real valued controlled stochastic differential equation (SDE, for short): 
\begin{align}\label{system0}
X_t&=x+\int_0^tb(s, X_s, u_s)\, ds+\int_0^t\sigma(s, X_s, u_s)^{\trans}\, dW_s, \quad t\in [0,T] 
\end{align}
with the time-point-wise expected path constraint
\begin{equation}\label{path constraint0}
\be[f(t, X_t)]\ge 0, \quad t\in [0,T]. 
\end{equation}
The cost functional reads
\begin{equation*} 
J(u):=\be\left[\int_0^T\ell(t,X_t,u_t)\, dt+g(X_T)\right]. 
\end{equation*}

The study of Stochastic Maximum Principles (SMPs, for short) is traced back to Bismut~\cite{Bismut1973,Bismut1976}, who introduced the notion of backward stochastic differential equations (BSDEs) to formulate the adjoint process and the stochastic Riccati equation, and was subsequently developed by Kushner~\cite{Kushner1972} and Haussmann~\cite{Haussmann1986}. At the early stage, these SMPs concerned only the stochastic systems where the control domain is convex or the diffusion coefficient does not contain control variable, and the proof only involves the first-order expansion. Peng~\cite{Peng1990} saw a breakthrough by establishing the SMP for the general stochastic optimal control problem where the control domain does not need to be convex and the diffusion coefficient can contain the control variable, where the second-order expansion and second-order backward stochastic differential equation (BSDE) are introduced. An extensive account of the progress on SMPs is available in Yong and Zhou~\cite{YongZhou}. 
Recently, SMP has found wide applications in probabilistic analysis of mean field games, and see the monograph~\cite{CarmonaDelarue2018}.

Our optimal stochastic control is featured by the inclusion of the expected path constraint. Our first aim is to establish a necessary condition (i.e. SMP) for this type of stochastic control problem, where the adjoint equation is a mean-reflected BSDE with the reflection being the consequence of the expected system path constraint. We note that a similar SMP has been already established by Frankowska et al.~\cite{Frankowska2019}; in contrast to theirs, our compensated
process $\mu_t$ is deterministic, which carries more information on the optimal control. Related results on optimal control of ordinary differential equations are referred to Dmitruk and Osmolovskii~\cite{DmitrukOsmolovskii2014} and Bourdin~\cite{Bourdin2016}

While applying this SMP to a stochastic control problem, a new type of
coupled reflected forward-backward stochastic differential equation
(FBSDE) appears:
\begin{equation*}
\begin{cases}
dX_t=(A_tX_t-B^{\trans}_tR^{-1}_t(B_tY_t+D^{\trans}_tZ_t))\, dt\\
\qquad\qquad+(C_tX_t-D_tR_t^{-1}(B_tY_t+D^{\trans}_tZ_t))^{\trans}\, dW_t, \\
dY_t=-(Q_tX_t+A_tY_t+C_t^{\trans}Z_t)\, dt+d\mu_{t}+Z^{\trans}_t\, dW_t, \\
\be[X_t]\ge L_{t}, \quad \int_0^T\left(\be[X_t]-L_{t}\right)\, d\mu_{t}=0, \\
X_{0}=x, \quad Y_T=GX_T, \quad \mu_{T}=0.
\end{cases}
\end{equation*}

This type of equation can be considered as an FBSDE counterpart of
BSDEs with mean reflection introduced by Briand et al. 
in \cite{BEH2018} and further studied by \cite{B2020}.
We will give a verification theorem and some well-solvability result concerning this new type of FBSDE.

The paper is organized as follows: after introducing some notation in the next subsection, we give the formulation of the problem in Section 2.
In Section 3, we apply Ekeland's variational principle to deduce the
stochastic maximum principle for the stochastic control problem.
In Section 4, we introduce the reflected FBSDE and show the verification theorem. The last two sections are devoted to the proof 
of uniqueness (Section 5) and of existence (Section 6).

\subsection{Notation}

Let $(W_t)_{0 \le t \le T}=(W_t^1,\cdots,W_t^m)_{0 \le t \le T}$
be an $m$-dimensional Brownian motion on a probability
space $(\Omega, \cF, \p)$. Denote by $\{\cF_t\}_{t\in [0,T]}$ the
augmented filtration generated by $(W_t)$.
Let $\R^+$ and $\R^-$, respectively, denote the sets of nonnegative and nonpositive real numbers. 
We write $x_+=\max\{x,\; 0\}$ and $x_-=\max\{-x,\; 0\}$ for $x\in\R$, the set of real numbers. 

We will often use vectors and matrices in this paper, where all vectors are column vectors. For a vector or matrix $M$,
denote by $M^{\trans}$ the transpose of $M$, and by $|M|=\sqrt{\sum_{i,j}m_{ij}^2}$ the Frobenius norm.

We will use the following notation. 
\begin{itemize}
\item $U$: a given closed convex subset of $\R^{l}$.

\item $\setu$: the set of $\{\cF_t\}_{t\in [0,T]}$-adapted functions $u(\cdot): [0,T]\times \Omega\to U$ such that $\E{ \int_0^T|u_t|^2\, dt} < \infty$.




\item $L^p_\cF([0, \, T]; \, \R^{k})$: the set of $\{\cF_t\}_{t\in [0,T]}$-adapted processes $f=(f_t^1,\cdots,f_t^k)_{0 \le t \le T}$ with $\E{ \int_0^T|f_t|^p\, dt} < \infty$. 


\item $L^\infty_\cF([0, \, T]; \, \R^{k})$: the set of essentially bounded $\{\cF_t\}_{t\in [0,T]}$-adapted $\R^{k}$-valued processes on $[0,T]$.
\item $L^\infty ([0, \, T]; \, \R^{k})$: the set of essentially bounded deterministic measurable $\R^{k}$-valued functions on $[0,T]$.

\item $C_\mathcal{F}([0,T];\,\R^{k})$: the Banach space of all continuous $\{\cF_t\}_{t\in [0,T]}$-adapted $\R^{k}$-valued processes $f$ on $[0,T]$ with a finite squared norm $\E{\max_{t\in [0,T]}|f_t|^2}$.
\item $\setmum$: the set of all nonpositive Radon measures on $[0,T]$. 
\item $\setmu$: the set of all nonnegative Radon measures on $[0,T]$.

\end{itemize}

For $\mu\in \setmu\cup\setmum$, we write
$$\mu_t=\mu([0,t])-\mu([0,T]).$$
Then the map $t\mapsto \mu_t$ is a c\`adl\`ag function on $[0,T]$
with $\mu_T=0$.

\section{Problem formulation}

Consider the following $\R$-valued controlled SDE: 
\begin{align}\label{system}
X_t&=x+\int_0^tb(s, X_s, u_s)\, ds+\int_0^t\sigma(s, X_s, u_s)^{\trans}\, dW_s, \quad t\in [0,T] 
\end{align}
with the time-point-wise expected path constraint
\begin{equation}\label{path constraint}
\be[f(t, X_t)]\ge 0, \quad t\in [0,T]. 
\end{equation}
The cost functional reads
\begin{equation*} 
J(u):=\be\left[\int_0^T\ell(t,X_t,u_t)\, dt+g(X_T)\right]. 
\end{equation*} 
In the above, $(b, \sigma): [0,T]\times \R\times U\times \Omega\to \R\times \R^{m}$, $f: [0,T]\times \R\to \R$, $\ell: [0,T]\times \R\times U\times \Omega\to \R$, and $g: \R\times \Omega\to \R$. 
\par
Let us assume the following conditions: 
\begin{description}
\item[(H1)] The maps $b$, $\sigma$, $f$, $\ell$ and $g$ are measurable. They are all continuously differentiable w.r.t. $(x,v)$. 
\item[(H2)] There exists a constant $L>0$ such that 
\begin{align*}
\begin{cases}
|b_{x}(t,x,v)|+|b_{v}(t,x,v)|+|\sigma_{x}(t,x,v)|+|\sigma_{v}(t,x,v)|\leq L,\\
|\ell_x(t,x,v)|+|\ell_v(t,x,v)|\leq L(1+|x|+|v|),\\
|g_x(x)|+|f_x(t,x)|\leq L (1+|x|),\\
|b(t,0,0)|+|\sigma(t,0,0)|+|\ell(t,0,0)|+|g(0)|+|f(t,0)|\leq L,
\end{cases}
\end{align*}
for any $(t,x, v,\omega)\in [0,T]\times\R\times U\times \Omega$.
\end{description}
We call a control $u\in\setu$ admissible if the SDE \eqref{system} admits a unique strong solution $X(\cdot)$ such that the constraint \eqref{path constraint} is satisfied.
The set of all admissible controls is denoted by $\setad$. 
We study the following optimal stochastic control problem 
\begin{equation}\label{target}
\min_{u\in\setad} J(u). 
\end{equation}

\section{Ekeland's variational principle and stochastic maximum principle}
We use Ekeland's variational principle to study the optimization problem~\eqref{target}. 
Before proceeding, we first present two technical lemmas. Denote by $C([0,T]; I)$ the set of all continuous functions $f: [0,T]\to I$ with $I=\R, \R^+$. Set $\setk:=C([0,T];\R^+) $, and define the distance function
$$
d_{\setk}(X):=\inf_{Y\in\setk} \|X-Y\|_{\infty}, \quad X\in C([0,T];\R),$$
with $\|\cdot\|_{\infty}$ being the maximal norm in  $C([0,T];\R)$. 

\begin{Lemma} \label{lem1}
For any $X\in C([0,T];\R)$, we have
$$
d_{\setk}(X)=\max_{t\in [0,T]} X_-(t) 
$$
\end{Lemma}

\begin{proof} First, since $X_{+}\in\setk$, we have 
$$d_{\setk}(X)=\inf_{Y\in\setk} \|X-Y\|_{\infty}\le \|X-X_{+}\|_{\infty}=\|X_-\|_{\infty}=\max_{t\in [0,T]} X_-(t).$$
It only remains to show the reverse inequality $d_{\setk}(X)\ge \|X_-\|_{\infty}.$ If $\|X_-\|_{\infty}=0$, then this inequality holds trivially. Otherwise, there is $t^*\in [0,T]$ such that \[X(t^*)=\min_{t\in [0,T]} X(t)=-\|X_-\|_{\infty}<0.\] 
So we have
$$
d_{\setk}(X)\ge \inf_{Y\in\setk} |X(t^*)-Y(t^*)|\ge |X(t^*)|=\|X_-\|_{\infty}.
$$
\end{proof} 

The subdifferential of the function $d_{\setk}$ at $X$, denoted by $\partial d_{\setk}(X)$, is defined to be the set of 
$\R$-valued Radon measures $K$ on $[0,T]$ such that
$$\langle K, f\rangle :=\int_{[0,T]} f(t) K(dt)\leq d_{\setk}(X+ f)-d_{\setk}(X), \quad \forall f\in \setk.$$

\begin{Lemma}\label{lem2}
For any $X\in C([0,T];\R)$, 
the set $\partial d_{\setk}(X)$ is not empty and 
$$\partial d_{\setk}(X)\subseteq \setmum$$ with $$\mbox{\rm supp } \partial d_{\setk}(X)\subseteq \mbox{\rm argmin } X. $$ Furthermore, if $X\not \in \setk$, we have $|K([0,T])|=1$ for any $K\in \partial d_{\setk}(X).$
\end{Lemma}

\begin{proof} 
We first show $\partial d_{\setk}(X)$ is not empty. If $d_{\setk}(X)=0$, then trivially $K\equiv 0\in\partial d_{\setk}(X)$. 
Otherwise $d_{\setk}(X)=-X(t^*)>0$ for some $t^*\in[0,T]$. 
Let $-K$ be the Dirac measure at $t^*$. 
Then by Lemma~\ref{lem1}, for any $h\in\setk$, 
\begin{align*}
d_{\setk}(X+h)-d_{\setk}(X)&=\max_{t\in [0,T]}(X(t)+h(t))_{-}+X(t^*)\\
&\geq \max_{t\in [0,T]}(-(X(t)+h(t)))+X(t^*)\geq -h(t^*)=\langle K, h\rangle. 
\end{align*}
Therefore, $K\in\partial d_{\setk}(X)$ and hence $\partial d_{\setk}(X)$ is not empty.

For any $(K,h)\in \partial d_{\setk}(X)\times \setk$, by Lemma~\ref{lem1},
$$
\langle K, h\rangle \le \lim_{\alpha\downarrow 0} { d_{\setk}(X+\alpha h)-d_{\setk}(X) \over \alpha}=\lim_{\alpha\downarrow 0}{ \|(X+\alpha h)_-\|_{\infty}-\|X_-\|_{\infty} \over \alpha}\le 0.
$$
Hence $K$ is non-positive. 
\par
If $t_0$ is not a minimum point of $X$, then 
$X$ has no minimum point on $[t_0-\varepsilon, t_0+\varepsilon]$ for sufficiently small $\varepsilon>0$. For any $h\in C([0,T];\R)$ with $\mbox{\rm supp}\;h\subset (t_0-\varepsilon/2, t_0+\varepsilon/2)$, we have $t_{0}$ is not a minimum point of $X\pm\alpha h$ for sufficiently small $\alpha>0$, so 
\[\|(X\pm\alpha h)_-\|_{\infty}=\|X_-\|_{\infty},\] which by definition implies $ \langle K, h\rangle=0$. Hence $\mbox{\rm supp } \partial d_{\setk}(X)\subseteq \mbox{\rm argmin } X$.

The last assertion is referred to \cite[Proposition 3.11, p. 146]{LiYong}. 
\end{proof} 

Let us first recall Ekeland's variational principle ( see~\cite[Theorem 1.1]{EK}).
\begin{Lemma}[Ekeland's variational principle]\label{evp}
Let $(V, d(\cdot,\cdot))$ be a complete metric space and $F(\cdot):V\to \mathbb R$ be a lower semi-continuous function, bounded from below. Suppose there exist $u\in V$ and $\varepsilon>0$ such that
$$F(u)\le \inf_{v\in V} F(v)+\varepsilon.$$
Then there exists $u_\varepsilon\in V$ such that
\begin{enumerate}
\item[(i)] $ F(u_\varepsilon)\le F(u)$,
\item[(ii)] $d(u,u_\varepsilon)\le \sqrt{\varepsilon}, \quad $ and 
\item[(iii)] $F(v)+\sqrt{\varepsilon}d(v,u_\varepsilon)\ge F(u_\varepsilon)$ for all $v\in V$.
\end{enumerate}
\end{Lemma}

We will work on the space $\setu$. To apply Ekeland's variational principle, we need to define a metric $d$ such that $(\setu,d)$ is a complete metric space. For this, set 
$$d(v,u)=\left( \be\left[\int_0^T|v(t)-u(t)|^2dt\right]\right)^{1/2}.$$
Then $(\setu, d(\cdot,\cdot))$ forms a complete metric space.

\par
Let $u^*\in\setad$ be an optimal control for problem~\eqref{target}. 
For $\varepsilon>0$ and $u\in\setu$, define the functional
\begin{align*}
J_\varepsilon(u)&=\left(\Big([J(u)-J(u^*)+\varepsilon]_+\Big)^2+\Big(\max_{t\in [0,T]}\big(\be[f(t,X_t^u)]\big)_-\Big)^2\right)^{1\over2}\\
&=\left(\Big([J(u)-J(u^*)+\varepsilon]_+\Big)^2+d_{\setk}^2\Big(\be[f(\cdot, X^u_\cdot)]\Big)\right)^{1\over2},
\end{align*} 
where the second equation is due to Lemma~\ref{lem1}. If $J_\varepsilon(u^\varepsilon)=0$, then $[J(u^\varepsilon)-J(u^*)+\varepsilon]_+=0$ and $u^\varepsilon\in\setad$, contradicting the optimality of $u^*$ to problem~\eqref{target}. So we have $J_\varepsilon(u^\varepsilon)> 0$.

Since 
$$ 
J_\varepsilon(u^*)=\varepsilon\le \inf_{u\in\setu} J_\varepsilon(u)+\varepsilon, 
$$
by Ekeland's variational principle Lemma \ref{evp}, we see that there is $u^\varepsilon\in\setu$ such that 
\begin{enumerate}
\item[(i)] $J_\varepsilon(u^\varepsilon)\le J_\varepsilon(u^*)$,
\item[(ii)] $d(u^*,u^\varepsilon)\le \sqrt{\varepsilon}, \quad $ and 
\item[(iii)] $J_\varepsilon(v)+\sqrt{\varepsilon}d(v,u^\varepsilon)\ge J_\varepsilon(u^\varepsilon)$ for all $v\in \setu$.
\end{enumerate}
The last assertion reads 
\begin{align}\label{Ekelandineq}
J_\varepsilon(u^\varepsilon)
=\min_{v\in \setu}(J_\varepsilon(v)+\sqrt{\varepsilon}d(v,u^\varepsilon)).
\end{align}

\par Let us establish the necessary condition for the optimization problem~\eqref{Ekelandineq}.
For any $v\in\setu$, $0<\alpha<1$ and $0<\varepsilon<1$, define 
$$
u^{\varepsilon,\alpha}=\alpha v+(1-\alpha)u^\varepsilon. 
$$
Denote by $X^{\varepsilon}$ and $X^{\varepsilon,\alpha}$ the trajectories corresponding to the controls $u^\varepsilon$ and $u^{\varepsilon,\alpha}$, respectively. By the Taylor expansion, we can identify $\delta X^\varepsilon$ and $\delta J(u^\varepsilon)$, which are independent of $\alpha$, such that for each fixed $\varepsilon$, 
$$
X^{\varepsilon,\alpha}=X^{\varepsilon}+\alpha\delta X^\varepsilon+o(\alpha)
$$
and 
$$
\Big([J(u^{\varepsilon,\alpha})-J(u^*)+\varepsilon]_+\Big)^2=\Big([J(u^\varepsilon)-J(u^*)+\varepsilon]_+\Big)^2+2\alpha[J(u^\varepsilon)-J(u^*)+\varepsilon]_+\delta J(u^\varepsilon)+o(\alpha)
$$
as $\alpha\to 0$. 
Then by \eqref{Ekelandineq}
\begin{align}
{(J_\varepsilon(u^{\varepsilon, \alpha}))^2-(J_\varepsilon(u^{\varepsilon}))^2\over J_\varepsilon(u^{\varepsilon, \alpha})+J_\varepsilon(u^{\varepsilon})}=J_\varepsilon(u^{\varepsilon, \alpha})-J_\varepsilon(u^{\varepsilon})\ge -\sqrt{\varepsilon} d(u^{\varepsilon, \alpha}, u^{\varepsilon})\ge -C\alpha \sqrt{\varepsilon}, 
\end{align}
where $C=d(v,u^{*})+1$ and the last inequality is due to 
\[\frac{1}{\alpha}d(u^{\varepsilon, \alpha}, u^{\varepsilon})=d(v,u^{\varepsilon})\leq d(v,u^{*})+d(u^*,u^\varepsilon)\le d(v,u^{*})+\sqrt{\varepsilon}\le d(v,u^{*})+1=C.\]
Therefore, 
$$
{(J_\varepsilon(u^{\varepsilon, \alpha}))^2-(J_\varepsilon(u^{\varepsilon}))^2\over \alpha \left(J_\varepsilon(u^{\varepsilon, \alpha})+J_\varepsilon(u^{\varepsilon})\right)}\ge -C \sqrt{\varepsilon}.
$$ 

Since $(C([0,T];\R),\|\cdot\|_{\infty}) $ is a separable Banach space, we know, see \cite{LiYong},
that there exists an equivalent norm, denoted by $\|\cdot\|_{0}$, such that the dual of 
$(C([0,T];\R), \|\cdot\|_{0})$ is strictly convex.  Any element $\mu$ of
$(C([0,T];\R), \|\cdot\|_{0})^*$ can still be identified with a  Radon measure on $[0,T]$. 
Since   $(C([0,T];\R), \|\cdot\|_{0})^*$ is strictly convex, $\partial d_{\setk}(X)$ is a singleton for any $X \notin \setk$. Furthermore, $d_{\setk}$ is G\^ateaux differentiable at any $X \notin \setk$. 
As $C$ does not depend on $\alpha$, letting $\alpha\downarrow 0$ in the last inequality, we obtain 
\begin{align}\label{ineq1}
{ [J(u^\varepsilon)-J(u^*)+\varepsilon]_+\delta J(u^\varepsilon)+d_{\setk}(\be[f(\cdot, X^\varepsilon_\cdot)])\int_0^T\be[f_{x}(t, X^\varepsilon_t)\delta X_t^\varepsilon]\, K^\varepsilon(dt)\over J_\varepsilon(u^{\varepsilon})}\ge -C \sqrt{\varepsilon}
\end{align}
where 
$K^\varepsilon\in \partial d_{\setk}(\be[f(\cdot, X^\varepsilon_\cdot)])\subseteq \setmum$.
Define 
$$
\lambda^\varepsilon:={ [J(u^\varepsilon)-J(u^*)+\varepsilon]_+\over J_\varepsilon(u^{\varepsilon})}\geq 0, \quad 
\mu^\varepsilon_t:=-{d_{\setk}(\be[f(\cdot, X^\varepsilon_\cdot)])K^\varepsilon([0,t])\over J_\varepsilon(u^{\varepsilon})}.
$$
By Lemma~\ref{lem2}, $|K^{\varepsilon}([0,T])|=1$ if $\be[f(\cdot, X^\varepsilon_\cdot)]\not \in \setk$, and $\mu^\varepsilon\equiv 0$ otherwise. Therefore, we have via a simple calculation
$$
|\lambda^\varepsilon|^2+|\mu^\varepsilon_T|^2=1.
$$
So there is a subsequence $\varepsilon_n\downarrow 0$ such that 
$$\lambda^{\varepsilon_n}\to \lambda\geq 0 \quad \text{and}\quad\mu^{\varepsilon_n}\to \mu,\quad \text{\rm $\star$-weakly in } C^*([0,T]; \R).$$
Since $\setk$ is obviously of finite-dimensional co-dimension in $C([0,T]; \R)$, in view of Lemma~\ref{lem2}, we have
$$
\lambda^\varepsilon\cdot 0+\langle \mu^\varepsilon, f\rangle \ge 0, \quad \forall f\in \setk.
$$
In view of \cite[Lemma 3.6, p. 142]{LiYong}, we have $(\lambda, \mu)\neq(0,0)$. 
\par
Set
\begin{align*}
b_x^\varepsilon(s)&:=b_x(s,X_s^\varepsilon, u^\varepsilon_s), & b_v^\varepsilon(s)&:=b_v(s,X_s^\varepsilon, u^\varepsilon_s),\\
b_x^*(s)&:=b_x(s,X_s^*, u^*_s), & b_v^*(s)&:=b_v(s,X_s^*, u^*_s), \\
\sigma_x^\varepsilon(s) &:=\sigma_x(s,X_s^\varepsilon, u^\varepsilon_s), & \sigma_v^\varepsilon(s)&:=\sigma_v(s,X_s^\varepsilon, u^\varepsilon_s),\\
\sigma_x^*(s)&:=\sigma_x(s,X_s^*, u^*_s),& \sigma_v^*(s)&:=\sigma_v(s,X_s^*, u^*_s),\\
\ell_x^*(s)&:=\ell_x(s,X_s^*, u^*_s), & \ell_v^*(s)&:=\ell_v(s,X_s^*, u^*_s),\\
\delta u^\varepsilon_s &:=v_s-u^\varepsilon_s, & \delta u^*_s &:=v_s- u^*_s.\\
\end{align*}
Then
\begin{align*}
\delta X_t^\varepsilon&=\int_0^t(b_x^\varepsilon(s)\delta X_s^\varepsilon+b_v^\varepsilon(s)\delta u^\varepsilon_s)\, ds+\int_0^t (\sigma_x^\varepsilon(s)\delta X_s^\varepsilon+\sigma_v^\varepsilon(s)\delta u^\varepsilon_s)^{\trans}\, dW_s,\\
\delta J(u^\varepsilon)&=\be\left[ \int_0^T\left(\ell_x^\varepsilon(s)\delta X_s^\varepsilon+\ell_v^\varepsilon (s)\delta u^\varepsilon_s\right) ds\right]
+\be[g_x(X_T^\varepsilon)\delta X_T^\varepsilon],\\
\delta X_t^*&=\int_0^t(b_x^*(s)\delta X_s^*+b_v^*(s)\delta u^*_s)\, ds+\int_0^t (\sigma_x^*(s)\delta X_s^*+\sigma_v^*(s)\delta u^*_s)^{\trans}\, dW_s,
\end{align*}
and
$$
\delta J(u^*)=\be\left[ \int_0^T(\ell_x^*(s)\delta X_s^*+\ell_v^* (s)\delta u^*_s)ds\right]
+\be[g_x(X_T^*)\delta X_T^*].
$$

By \eqref{ineq1},
$$
\lambda^\varepsilon \delta J(u^\varepsilon)+\int_0^T\be[f_{x}(t, X^\varepsilon_t)\delta X_t^\varepsilon]\, d\mu_t^\varepsilon\ge -C \sqrt{\varepsilon}.
$$
As $C$ does not depend on $\varepsilon$, letting $\varepsilon\downarrow 0$, 
\begin{align}\label{ineq2}
\lambda \delta J(u^*)+\int_0^T\be[f_{x}(t, X^*_t)\delta X_t^{*}]\, d\mu_t\ge 0.
\end{align}
Denote by $(Y, Z)$ the unique solution of the following BSDE
\begin{eqnarray}\label{BSDE}
Y_t=&\displaystyle \lambda g_x(X_T^*)+\int_t^T(b_x^{*}(s)Y_s+\sigma_x^{*}(s)^{\trans}Z_s+\lambda \ell_x^{*}(s))\, ds \\
&\displaystyle +\int_t^T f_x(s, X^*_s) d\mu_{s}-\int_t^TZ_s^{\trans}\, dW_s. \nonumber
\end{eqnarray}
We have the following stochastic maximum principle.

\begin{Theorem} Let $u^*\in\setad$ be an optimal control for problem~\eqref{target}. Then, there is $(\lambda, \mu)\in [0,1]\times \setmu$ such that (i) $(\lambda, \mu)\neq(0,0)$ and (ii)  the following maximum condition is satisfied:
$$
\min_{v\in U}\; \left\{\langle Y_t, b_v^*(t)(v-u^*_t)\rangle+\langle Z_t, \sigma_v(t) (v-u^*_t)\rangle+\lambda \ell_v^*(t)(v-u^*_t)\right\}=0, \quad \mathrm{a.e.}\; t\in[0,T], 
$$
where the pair $(Y, Z)$ is the unique solution of  BSDE~\eqref{BSDE}. 
\end{Theorem}

\begin{proof}
By \eqref{ineq2}, we have
\begin{align*} 
0&\le\be\left[ \int_0^T\lambda(\ell_x^*(s)\delta X_s^*+\ell_v^* (s)\delta u^*_s)ds+\lambda g_x(X_T^*)\delta X_T^* +\int_0^Tf_{x}(t, X^*_t)\delta X_t^{*}\, d\mu_t\right]\\
&=\be\left[\int_0^T\lambda \ell_v^* (s)\delta u^*_sds-\int_0^T \delta X_s^*d Y_s+Y_T\delta X_T^*-\int_0^T (b_x^{*}(s) Y_s+\sigma_x^{*}(s)^{\trans}Z_s) \delta X_s^*ds\right]\\
&=\be\bigg[\int_0^T\lambda \ell_v^* (s)\delta u^*_sds+\int_0^T Y_s (b_x^*(s)\delta X_s^*+b_v^*(s)\delta u^*_s)\, ds\\
&\qquad\;\;+ \int_0^TZ_s^{\trans}(\sigma_x^*(s)\delta X_s^*+\sigma_v^*(s)\delta u^*_s)\, ds-\int_0^T(b_x^{*}(s) Y_s+\sigma_x^{*}(s)^{\trans}Z_s) \delta X_s^*ds\bigg]\\
&=\be\bigg[\int_0^T\lambda \ell_v^* (s)\delta u^*_s+Y_s b_v^*(s)\delta u^*_s +Z_s^{\trans}\sigma_v^*(s)\delta u^*_sds\bigg]. 
\end{align*}
This implies the desired result. 
\end{proof}

\section{LQ stochastic control problem with expected path constraints}

We now study an LQ stochastic control problem with an expected path constraint. The dynamic of the state process is governed by the SDE
\begin{equation}\label{state}
dX_{t}=(A_{t}X_{t}+B _{t}^{\trans}u_{t}) dt+(C_{t}X_{t}+D_{t}u_{t})^{\trans} dW _{t}. 
\end{equation}
Here the state process $X$ is one-dimensional and the control $u\in\setu$ is $l$-dimensional.
The coefficient matrices $A$, $B$, $C$, $D$ are essentially bounded adapted processes of proper sizes. 

Let $\setad$ be the set of all controls $u\in \setu$ such that the pair $(u,X)$ solves equation~\eqref{state} with the initial value $X(0)=x$, and satisfies the following expected path constraint 
\begin{equation}\label{constraint}
\be[X_t]\ge L_{t},\quad \forall\; t \in [0,T].
\end{equation}
Here, $L$ is a given deterministic continuous function. 
Introduce the constrained problem
\begin{equation}\label{cost}
\min_{u\in\setad} J(u):={1\over 2} \be\left[\int_0^T (Q_{t}X^2_t+u_t^{\trans} R_{t}u_t)\, dt+GX^2_T\right],
\end{equation}
and denote by $V(x)$ its optimal value.

\begin{Assumption}\label{assump0}
We have $Q\ge 0, G\ge 0$, and $R \ge\delta I_l$ uniformly in $(t,\omega)$ for some $\delta>0$.
\end{Assumption}

To guarantee that the admissible set $\setad$ is not empty, we put the following assumption. 
\begin{Assumption}\label{assump1}
There exist a control $u^{a}\in\setu$ and a constant $\ep>0$ such that $(u^{a},X^{a})$, which solves equation~\eqref{state} with $X^{a}_{0}=x$, satisfies 
$\be\left[X_t^{a}\right]> L_{t}+\ep$ for all $t\in[0,T].$
\end{Assumption}

The last assumption holds true if $L_{t}<-\ep$ for all $ t\in [0,T]$ and $x\ge 0$. In fact, it suffices to choose $u^{a}=0$.

\begin{Remark}
Suppose there is $u^{a}\in\setu$ such that $(u^{a}, X^{a})$ solves equation~\eqref{state} with $X^{a}_{0}<x$ and satisfies $\be\left[X_t^{a}\right]\ge L_{t}$ for all $t\in[0,T].$ Then Assumption~\ref{assump1} holds. 
In fact, suppose $(u^{a}, X)$ solves equation~\eqref{state} with $X_{0}=x$. Then by the strict monotonicity of the SDE~\eqref{state} with respect to the initial value, we have $\be[X_t-X_t^{a}]>\ep$ for some $\ep>0$. So $\be\left[X_t\right]>L_{t}+\ep$ for all $t\in[0,T].$
\end{Remark}

We assume that Assumptions \ref{assump0} and \ref{assump1} hold in the rest of the paper. 

\subsection{Existence, uniqueness and approximation of the optimal control}\label{existence}

\begin{Lemma} Let Assumptions~\ref{assump0} and \ref{assump1} be satisfied. Then, Problem \eqref{cost} has a unique optimal control. 
\end{Lemma}

\begin{proof} We first show that problem \eqref{cost} has an optimal solution. In fact, from Assumption~\ref{assump1}, we see that there is a minimizing sequence $\{v^n, n=1,2,\cdots\}$ in the set $\setad$. It suffices to prove that $\{v^n, n=1,2,\cdots\}$ is a Cauchy sequence in the Banach space $L^2_\cF([0, \, T]; \, \R^{l})$, since its limit still lies in $\setad$. We have 
$$
\lim_{n\to \infty}J(v^n)=V(x), \quad v^{n,k}:=\frac{1}{2}(v^n+v^k)\in \setad,\quad X^{n,k}=\frac{1}{2}(X^n+X^k)
$$
where $X^n$ and $X^{n,k}$ are the state processes under the admissible controls $v^n$ and $v^{n,k}$, respectively. 
Therefore, $J(v^{n,k})\ge V(x)$, and the parallelogram rule holds:
\begin{align}\label{parallelogram}
&\quad\;\frac{1}{4} \be\left[\int_0^T [Q_{t}(X_t^n-X_t^k)^2+(v_t^n-v_t^k)^{\trans} R_{t}(v_t^n-v_t^k)]\, dt+G(X_T^n-X_T^k)^2\right]\nonumber\\
&=J(v^n)+J(v^k)-2J(v^{n,k}) \le J(v^n)+J(v^k)-2V(x). 
\end{align}
Hence, we have 
$$
\frac{1}{4}\delta \|v^n-v^k\|^2\le J(v^n)+J(v^k)-2V(x) \to 0,  \quad \text{ as $n,k\to \infty$}, 
$$
and then $\{v^n, n=1,2,\cdots\}$ is a Cauchy sequence in the Banach space $L^2_\cF([0, \, T]; \, \R^{l})$. 

Uniqueness of optimal control can be proved in a similar way via the parallelogram rule. 
\end{proof}

On the other hand, consider the following unconstrained problem for each $n>0$, 
\begin{equation}\label{unconstrainedproblem}
\min_{u\in\setu} J_n(u):=\frac12\be \left[ \int_0^TQ_{t}X^2_{t}+u_t^{\trans} R_{t}u_t\, dt+GX^2_T\right]+\frac12n\int_0^T\left[(\be[X_t]-L_{t})_-\right]^2dt,
\end{equation}
where the state process $X$ solves equation~\eqref{state}. Problem~\eqref{unconstrainedproblem} is a stochastic linear-convex optimal control problem, which admits a unique solution (see a similar proof of Yong and Zhou~\cite[Theorem 5.2, page
68]{YongZhou}).

Let $V_n(x)$ be the optimal value function of \eqref{unconstrainedproblem}. Then, for any control $u\in\setad$, we have
\begin{equation}
J_n(u)=J(u),
\end{equation}
which leads to
\begin{equation}\label{vnlev}
V_n(x) \le V (x). 
\end{equation}

\begin{Lemma} Let $(\bar u^{n}, \overline{X}^{n})$ be the optimal pair of the unconstrained problem~\eqref{unconstrainedproblem}. 
Then, $\bar u^{n}$ converges strongly to the optimal control of the constrained problem \eqref{cost}. 
\end{Lemma}

\begin{proof} We have 
\begin{align}\label{n costa}
J_n(\bar{u}^n) &=\frac12\be \left[\int_0^TQ_{t}(\overline{X}^n_t)^2+(\bar{u}^n_{t})^{\trans} R_{t}\bar{u}^n_{t}\, dt+G(\overline{X}^n_T )^2\right] \nonumber\\
&\quad\; +\frac12n \int_0^T[(\be[ \overline{X}^{n}_{t}]-L_{t})_-]^2dt\le J_n(u^{a})=J( u^{a}), 
\end{align}
where $u^a$ is given in Assumption~\ref{assump1}. As $R \ge \delta I_l$, it follows that the sequence $\bar{u}^n$ is bounded in $L^2_\cF(0, \, T; \, \R^l)$. Consequently, it has a subsequence (still denoted by $\bar{u}^n$) which weakly converges to some control $\bar{u}^\infty \in L^2_\cF(0, \, T; \, \R^l)$.

By Mazur's theorem, there exist real numbers $\epsilon_{n,k}\geq0$ such that $\sum_{k\geq 0} \epsilon_{n,k}=1$ for every $n$, and the sequence 
\[\bar{v}^n=\sum_{k\geq 0} \epsilon_{n,k} \bar{u}^{n+k}\]
strongly converges to $\bar{u}^\infty$. Let $X^{\bar{v}^n}$ and $\overline{X}^\infty$ denote, respectively, the trajectories under the controls $\bar{v}^n$ and $\bar{u}^\infty$. Then the sequence $X^{\bar{v}^n}$ converges to $\overline{X}^\infty$ strongly in $C_\mathcal{F}([0,T];\R)$. In particular 
\[\overline{X}^\infty_t=\lim_n X^{\bar{v}^n}_t=\lim_n \sum_{k\geq 0}\epsilon_{n,k} \overline{X}^{n+k}_t,\quad t\in[0,T].\]
Dividing both sides of \eqref{n costa} by $\frac{n}2$ and letting $n$ go to $\infty$, we deduce from the convexity of the map $x\mapsto (x_{-})^2$ and Fatou's lemma that
\begin{align*}
\int_0^T[(\be[\overline{X}^\infty_t]-L_{t})_{-}]^2 dt 
&\le \liminf_{n\to \infty} \sum_{k\geq 0}\epsilon_{n,k}\int_0^T[(\be[\overline{X}^{n+k}_t]-L_{t})_{-}]^2 dt\\
&\le \liminf_{n\to \infty} \sum_{k\geq 0}\epsilon_{n,k}\frac{2}{n}J( u^{a})= \liminf_{n\to \infty} \frac{2}{n}J( u^{a})=0.
\end{align*}
Since $\be[\overline{X}^\infty_{\cdot}]$ and $L_{\cdot}$ are continuous, we conclude $\be[\overline{X}^\infty_t]\ge L_{t}$ holds for all $t\in[0,T]$. This means $ (\bar{u}^\infty,\overline{X}^\infty)$ is an admissible pair for the constrained problem~\eqref{cost}, so
\begin{align} J(\bar{u}^\infty) \ge V (x).\label{existeq0}\end{align}
The convexity of the map $x\mapsto (x_{-})^2$ and Fatou's lemma also give
\begin{align}
&\limsup_{n\to \infty}\be \left[\int_0^TQ_{t}(\overline{X}^{n}_t)^2 \, dt+G(\overline{X}^{n}_T )^2\right]\nonumber\\
\geq & \liminf_{n\to \infty} \sum_{k\geq 0} \epsilon_{n,k} \be \left[\int_0^TQ_{t}(\overline{X}^{n+k}_t)^2 \, dt+G(\overline{X}^{n+k}_T )^2\right] \nonumber\\
\geq &\be \left[\int_0^TQ_{t}(\overline{X}^\infty_t)^2 \, dt+G(\overline{X}^\infty_T )^2\right].
\label{existeq1}
\end{align}
Thanks to the weak convergence of $\bar{u}^n$ to $\bar{u}^\infty$, 
\begin{align}
\liminf_{n\to \infty} \be \left[\int_0^T (\bar{u}^n_{t})^{\trans} R_{t}\bar{u}^n_{t}\, dt \right] \ge \be \left[\int_0^TQ_{t} (\bar{u}^\infty_{t})^{\trans} R_{t}\bar{u}^\infty_{t}\, dt \right]. \label{existeq2}
\end{align}
The above estimates yield 
\begin{align}
\lim_{n\to\infty}J_{n}(\bar{u}^n) &\geq \limsup_{n\to \infty} \be \left[\int_0^TQ_{t}(\overline{X}^n_t)^2+(\bar{u}^n_{t})^{\trans} R_{t}\bar{u}^n_{t}\, dt+G(\overline{X}^n_T )^2\right]\nonumber\\
&\ge \be \left[\int_0^TQ_{t}(\overline{X}^\infty_t)^2+(\bar{u}^\infty_{t})^{\trans} R_{t}\bar{u}^\infty_{t}\, dt+G(\overline{X}^\infty_T )^2\right]\nonumber \\
&=J(\bar{u}^\infty) \ge V (x).\label{existeq3}
\end{align}
But \eqref{vnlev} gives 
\[V(x) \ge V_n (x)=J_{n}(\bar{u}^n),\]
so all the inequalities in \eqref{existeq0}-\eqref{existeq3} are equations. In particular, $(\bar{u}^\infty, \overline{X}^\infty)$ is the optimal pair of the constrained problem~\eqref{cost} as \eqref{existeq0} is an equation. By the weak convergence and norm convergence \eqref{existeq2} we conclude $R^{\frac12}\bar{u}^n$ strongly converges to $R^{\frac12}\bar{u}^\infty$ in the space $L^2_{\cal F}([0,T];\mathbb R)$. As $R \ge \delta I_l$, $\bar{u}^n$ strongly converges to $\bar{u}^\infty$ in the space $L^2_{\cal F}([0,T];\mathbb R)$. Consequently, $\overline{X}^n$ strongly converges to $\overline{X}^\infty$ in the space $L^2_{\cal F}([0,T];\mathbb R)$.
As a byproduct of \eqref{existeq3}, we have 
\begin{align}\label{limV}
\lim_{n\to \infty}\,\, n\int_0^T[(\be[ \overline{X}^{n}_{t}]-L_{t})_-]^2dt=0. 
\end{align}

Finally we note that as the optimal control is unique, the whole sequence $\bar{u}^n$ strongly converges to $\bar{u}^\infty$ in the space $L^2_{\cal F}([0,T];\mathbb R)$.
\end{proof}

\subsection{Verification theorem}
In this section, we express the unique optimal control for problem \eqref{cost} with the solution of a reflected FBSDEs. 
\par
We say that $(X,Y,Z,\mu)\in L^2_\cF([0, \, T]; \, \R^3)\times \setmu$ is a solution of the following reflected FBSDEs
\begin{equation}\label{FBSDEs1}
\begin{cases}
dX_t=(A_tX_t-B^{\trans}_tR^{-1}_t(B_tY_t+D^{\trans}_tZ_t))\, dt\\
\qquad\qquad+(C_tX_t-D_tR_t^{-1}(B_tY_t+D^{\trans}_tZ_t))^{\trans}\, dW_t, \\
dY_t=-(Q_tX_t+A_tY_t+C_t^{\trans}Z_t)\, dt+d\mu_{t}+Z^{\trans}_t\, dW_t, \\
\be[X_t]\ge L_{t}, \quad \int_0^T\left(\be[X_t]-L_{t}\right)\, d\mu_{t}=0, \\
X_{0}=x, \quad Y_T=GX_T, \quad \mu_{T}=0,
\end{cases}
\end{equation} 
if it satisfies the above FBSDEs.

\begin{Theorem} \label{varify1}
Suppose that $(\overline{X},\overline{Y},\overline{Z},\overline{\mu})$ is a solution of the reflected FBSDEs \eqref{FBSDEs1}. 
Then \[\overline{u}:=-R^{-1}(B\overline{Y}+D^{\trans} \overline{Z})\] is an optimal control for problem~\eqref{cost}, and the optimal value is 
\begin{equation}
J(\overline{u})={1\over2} \overline{Y}_0x+{1\over2}\int_0^TL_{t}\,d\overline{\mu}_{t}. 
\end{equation}
\end{Theorem}

\begin{proof} 
Note that $(\overline{X}, \overline{u})$ solves the equation~\eqref{state} with the initial value $x$.  So for any $u\in\setad$, 
\begin{align}
&\quad\;J(u)-J(\overline{u})\nonumber\\
&={1\over2} \be[GX_T^2-G\overline{X}_T^2]+{1\over2} \be\int_0^T[Q_{t}X_t^2-Q_{t}\overline{X}_t^2]dt+{1\over2} \be\int_0^T[u_t^{\trans} R_{t}u_t-\overline{u}_t^{\trans} R_{t}\overline{u}_t]dt\nonumber \\
&=\be[G\overline{X}_T(X_T-\overline{X}_T)]+\be\int_0^TQ_{t}\overline{X}_t(X_t-\overline{X}_t)\, dt+\be\int_0^T\overline{u}_t^{\trans} R_{t}(u_t-\overline{u}_t)\, dt\nonumber\\
&\quad\;+{1\over2}\be[G(X_T-\overline{X}_T)^2]+{1\over2}\be\int_0^TQ_{t}(X_t-\overline{X}_t)^2\, dt+{1\over2}\be\int_0^T(u_t-\overline{u}_t)^{\trans} R_{t}(u_t-\overline{u}_t) \, dt.\nonumber\\
&\geq \be[G\overline{X}_T(X_T-\overline{X}_T)]+\be\int_0^TQ_{t}\overline{X}_t(X_t-\overline{X}_t)\, dt+\be\int_0^T\overline{u}_t^{\trans} R_{t}(u_t-\overline{u}_t)\, dt.\label{ineq3}
\end{align}
Applying It\^o's formula, we have
\begin{align*}
d(\overline{Y}_t(X_t-\overline{X}_t))&=-Q_{t}\overline{X}_t(X_t-\overline{X}_t)\,dt+(X_t-\overline{X}_t)\, d\overline{\mu}_{t}+\overline{Y}_tB^{\trans}_t(u_t-\overline{u}_t)\, dt \\
&\quad\;+\overline{Z}^{\trans}_t [D(u_t-\overline{u}_t)]\, dt+(X_t-\overline{X}_t)\overline{Z}^{\trans}_t dW_t \\
&\quad\;+\overline{Y}_t[C_t(X_t-\overline{X}_t)+D(u_t-\overline{u}_t)]^{\trans} dW_t. 
\end{align*}
Integrating both sides and taking the expectation (also noting that the local martingale is in fact a martingale (see Bismut~\cite[Proposition I-1, p. 387]{Bismut1973})), we have 
\begin{align*}
&\quad\;\be[G\overline{X}_T(X_T-\overline{X}_T)]+\be\int_0^T[Q_{t}\overline{X}_t(X_t-\overline{X}_t)]\, dt \\
&=\int_0^T(\be[X_t]-\be[\overline{X}_t])\, d\overline{\mu}_{t}+E\left[\int_0^T\langle B_{t}\overline{Y}_{t}+D^{\trans}_{t} \overline{Z}_{t}, u_t-\overline{u}_t\rangle dt\right].
\end{align*}
Thanks to $\overline{u}=-R^{-1}(B\overline{Y}+D^{\trans} \overline{Z})$ and \eqref{ineq3}, 
\begin{align*}
&\quad\; J(u)-J(\overline{u}) \\ 
&\ge \be[G\overline{X}_T(X_T-\overline{X}_T)]+\be\int_0^T[Q_{t}\overline{X}_t(X_t-\overline{X}_t)]\, dt+\be\left[\int_0^T\langle R_{t}\overline{u}_{t}, u_t-\overline{u}_t\rangle dt\right] \\
&=\int_0^T(\be[X_t]-L_{t})\, d\overline{\mu}_{t}-\int_0^T(\be[\overline{X}_t]-L_{t})\, d\overline{\mu}_{t}=\int_0^T(\be[X_t]-L_{t})\, d\overline{\mu}_{t}\ge 0,
\end{align*}
where the last inequality is due to the constraint \eqref{constraint} and $\overline{\mu}\in\setmu$. 

Again, using It\^o's formula, we have
\begin{align*}
d(\overline{Y}_t\overline{X}_t)&=-Q_{t}\overline{X}_t^{2}\,dt+\overline{X}_t\, d\overline{\mu}_{t}+\overline{Y}_tB^{\trans}_t\overline{u}_t\, dt+\overline{Z}^{\trans}_t D_{t}\overline{u}_t\, dt\nonumber\\
&\quad\; +\overline{X}_t\overline{Z}_t^{\trans} dW_t+\overline{Y}_t(C_t\overline{X}_t+D\overline{u}_t)^{\trans} dW_t.\nonumber
\end{align*}
Note that the local martingale is in fact a martingale (see Bismut~\cite[Proposition I-1, p. 387]{Bismut1973}). 
Therefore, integrating both sides yields
\begin{align*}
\be[G\overline{X}_T^2]-\overline{Y}_0x
&=-\be\left[\int_0^TQ_{t}\overline{X}_t^2\, dt \right] +\int_0^T\be[\overline{X}_t]\, d\overline{\mu}_{t}+\be\left[\int_0^T\langle B^{\trans}_{t} \overline{Y}_{t}+D^{\trans}_{t} \overline{Z}_{t}, \overline{u}_t\rangle dt\right] \\
&=-\be\left[\int_0^TQ_{t}\overline{X}_t^2\, dt \right]+\int_0^TL_{t}\, d\overline{\mu}_{t}-\be\left[\int_0^T\langle R_{t}\overline{u}_t, \overline{u}_t\rangle dt\right]. 
\end{align*}
Thus, we proved the desired expression for the optimal value $J(\overline{u})$.
\end{proof}

In the rest of the paper, we focus on solution of the reflected FBSDEs \eqref{FBSDEs1}. 
The main result is stated as follows. 
\begin{Theorem}\label{existanduniqueFBSDEs}
If $A$ is deterministic, $B^{\trans}B$ is invertible and
$(B^{\trans}B)^{-1}$ is bounded, 
then the reflected FBSDEs \eqref{FBSDEs1} admits a unique solution.
\end{Theorem}

\begin{proof}
This is an immediate consequence of Propositions~\ref{FBSDEunique} and ~\ref{existenceFBSDEs} in the subsequent sections. 
\end{proof}

We will prove the uniqueness and existence in the following two sections respectively.

\section{Uniqueness of the solution for the reflected FBSDEs \eqref{FBSDEs1}}

\begin{Proposition}\label{FBSDEunique} 
Let $(X, Y, Z, \mu)$ and $(\widehat{X}, \widehat{Y}, \widehat{Z}, \widehat{\mu})$ be two solutions for the reflected FBSDEs \eqref{FBSDEs1}. Then $X=\widehat{X}$ and $BY+D^{\trans} Z=B\widehat{Y}+D^{\trans} \widehat{Z}$. Furthermore, $(X, Y, Z, \mu)=(\widehat{X}, \widehat{Y}, \widehat{Z}, \widehat{\mu})$ if $A$ is deterministic and $\be[{B^{\trans}B}]>0$. 
\end{Proposition}

\begin{proof} We denote by $(\widetilde{X}, \widetilde{Y}, \widetilde{Z}, \widetilde{\mu})$ the difference of $(X, Y, Z, \mu)$ and $(\widehat{X}, \widehat{Y}, \widehat{Z}, \widehat{\mu})$. Set
$$
\overline{u}:=-R^{-1}(BY+D^{\trans} Z), \quad 	\widehat{u}:=-R^{-1}(B\widehat{Y}+D^{\trans} \widehat{Z}),
\quad \widetilde{u}:=\overline{u}-\widehat{u}. 
$$

We now show the first assertion, that is, $\widetilde{X}=0$ and $\widetilde{u}=0$. 
By \eqref{FBSDEs1}, 
\begin{align}\label{sde1}
d\widetilde{X}&=(A \widetilde{X}+B^{\trans} \widetilde{u})\, dt+(C \widetilde{X}+D \widetilde{u})^{\trans} dW_t, \\
d\widetilde{Y}&=- (Q \widetilde{X}+A \widetilde{Y}+C^{\trans}\widetilde{Z})\, dt+d\widetilde{\mu}_{t}+\widetilde{Z}_t^{\trans} dW_t.\nonumber
\end{align}
Using It\^o's formula, we have
\begin{align*}
d(\widetilde{X}_t\widetilde{Y}_t)&=B^{\trans} \widetilde{u}\widetilde{Y}\, dt+(C\widetilde{X}+D\widetilde{u})^{\trans} \widetilde{Y} dW_t\nonumber\\
&\quad\;-\widetilde{X}(Q\widetilde{X}\, dt-d\widetilde{\mu}_{t}-\widetilde{Z}^{\trans} dW_t)+\widetilde{Z}_t^{\trans} D\widetilde{u} dt.
\end{align*}
Integrating both sides and taking the expectation, since the local martingale is in fact a martingale (see Bismut~\cite[Proposition I-1, p. 387]{Bismut1973}), we have the duality formula 
\begin{align*}
&\quad\;\be[G\widetilde{X}_T^2]+\be\int_0^T Q \widetilde{X}^2_t dt \nonumber\\
&=\be\int_0^T\widetilde{X} d\widetilde{\mu}_{t}+\be\int_0^T\langle B \widetilde{Y}+D^{\trans} \widetilde{Z}, \widetilde{u}\rangle dt\nonumber\\
&=\int_0^T(\be[X_{t}]-L_{t})d\widetilde{\mu}_{t}-\int_0^T(\be[\widehat{X}_{t} ]-L_{t})d\widetilde{\mu}_{t}-\be\int_0^T\langle R \widetilde{u}_{t}, \widetilde{u}_{t}\rangle dt.\nonumber
\end{align*}
Hence,
\begin{align*}
&\quad\;\be[G\widetilde{X}_T^2]+\be\int_0^T Q \widetilde{X}^2_t dt+\be\int_0^T\langle R \widetilde{u}, \widetilde{u}\rangle dt\nonumber\\
&=-\int_0^T(\be[X_{t}]-L_{t})d\widehat{\mu}_{t}-\int_0^T(\be[\widehat{X_{t}} ]-L_{t})d\mu_{t}\le 0.\nonumber
\end{align*}
Because $G\ge 0$, $Q\ge 0$, $R \ge \delta I_l$, it follows $\widetilde{u}=0$. 
Consequently, \eqref{sde1} reduces to 
$d\widetilde{X}=A \widetilde{X}\,dt+ (C \widetilde{X})^{\trans} dW_t$. 
Together with $\widetilde{X}_0=0$ we infer that $\widetilde{X}=0$. 
This completes the proof of the first assertion. 

\par
Now suppose $A$ is deterministic. 
Let $U_{t}=e^{\int_0^t A_rdr}>0$. Then $U$ is deterministic and $d U=AUdt$. By \eqref{sde1}, 
\begin{align*} 
d (U\widetilde{Y})&=- UC^{\trans}\widetilde{Z}\, dt+Ud\widetilde{\mu}_{t}+U \widetilde{Z}^{\trans} dW_t=U\widetilde{Z}^{\trans} d \widetilde{W}_t+U d\widetilde{\mu}_{t}, \label{tildeyeq}
\end{align*}
where $\widetilde{W}_{t}=W_t-\int_0^t C_{s} ds$ is a Brownian motion under some probability measure $\widetilde{\p}\sim\p$. 
This means 
\begin{equation*} 
U_{t}\widetilde{Y}_{t}-\int_{0}^{t}U_{s}d\widetilde{\mu}_{s}=U_{0}\widetilde{Y}_0+\int_0^tU_s\widetilde{Z}_s^{\trans} d \widetilde{W}_s
\end{equation*}
is a martingale under $\widetilde{\p}$. But the value of this martingale at $t=T$ is 
\[U_{T}\widetilde{Y}_{T}-\int_{0}^{T}U_{s}d\widetilde{\mu}_{s}=
U_{T}G\widetilde{X}_{T}-\int_{0}^{T}U_{s}d\widetilde{\mu}_{s}
=-\int_{0}^{T}U_{s}d\widetilde{\mu}_{s},\]
a constant, so it is a constant martingale. 
Hence $\widetilde{Z}=0$ and consequently, 
\begin{equation}\label{martingale1}
\widetilde{Y}_{t}=U_{t}^{-1}\left(U_0\widetilde{Y}_0+\int_{0}^{t}U_{s}d\widetilde{\mu}_{s}\right) 
\end{equation}
is a deterministic function. 
From $\widetilde{u}=0$ and $\widetilde{Z}=0$, 
we get $B\widetilde{Y}=0$. Thus 
$$0=\be[B^{\trans}B\widetilde{Y}]=\be[{B^{\trans}B}]\widetilde{Y}. $$
If $\be[{B^{\trans}B}]>0$, then $\widetilde{Y}=0$, and consequently by \eqref{martingale1}, $\widetilde{\mu}=0$. The second assertion is thus proved.
\end{proof} 

\begin{Remark}
If the last condition in the above theorem does not hold, then the uniqueness can fail. 
For instance, when $A$ is deterministic and $B=0$, we may get infinite many solutions $(Y+(k-1)\widetilde{Y},\; k{\mu})$ from a solution $(Y,\mu)$ by setting {$k>0$} and 
\[\widetilde{Y}_t=-\int_{t}^{T}e^{\int_t^s A_rdr}d\mu_s.\]
\end{Remark}

\section{Existence of the solution for the reflected FBSDEs \eqref{FBSDEs1}}

\begin{Proposition}\label{existenceFBSDEs}
If $B^{\trans}B$ is invertible and
$(B^{\trans}B)^{-1}$ is bounded.
Then the reflected FBSDEs \eqref{FBSDEs1} has a solution.
\end{Proposition}
We use the penalization method to prove the existence.
The proof is given in the subsequent two subsections. 

\subsection{Approximation}
For any $n \in\mathbb{N}$, consider the following penalized FBSDEs:
\begin{equation}\label{FBSDEs2}
\begin{cases}
dX^n=(AX^n-B^{\trans} R^{-1}(BY^n+D^{\trans} Z^n))\,dt\\
\qquad\qquad+(CX^n-DR^{-1}(BY^n+D^{\trans} Z^n))^{\trans} dW_t, \\
dY^n=-(QX^n+AY^n+C^{\trans} Z^n)\,dt+n(\be[X^n]-L)_-dt+(Z^n)^{\trans} dW_t, \\
X^n_{0}=x, \quad Y_T^n=GX_T^n;
\end{cases}
\end{equation}
It is a Mckean-Vlasov FBSDEs, and is actually the Hamiltonian system of the optimal control of Problem \eqref{unconstrainedproblem}. 
We call $(X^n, Y^n, Z^n) \in L^2_\cF([0, \, T]; \, \R^3)$ a solution to the FBSDEs \eqref{FBSDEs2} if it satisfies \eqref{FBSDEs2}.

\begin{Lemma}\label{penalexist}
The penalized FBSDEs \eqref{FBSDEs2} admits a unique solution.
\end{Lemma}
\begin{proof} Since the optimal control problem~\eqref{unconstrainedproblem} has an optimal control $\bar u$, we have the existence of FBSDEs~\eqref{FBSDEs2} immediately from the stochastic maximum principle for optimally controlled Mckean-Vlasov SDEs. 

We now turn to the proof of uniqueness. We suppress the superscript $n$ here for simplicity.
We denote by $(\widetilde{X}, \widetilde{Y}, \widetilde{Z})$ the difference of two solutions $(X, Y, Z)$ and $(\widehat{X}, \widehat{Y}, \widehat{Z})$ to ~\eqref{FBSDEs2}.
Set
$$
\overline{u}:=-R^{-1}(BY+D^{\trans} Z), \quad 	\widehat{u}:=-R^{-1}(B\widehat{Y}+D^{\trans} \widehat{Z}),
\quad \widetilde{u}:=\overline{u}-\widehat{u}. 
$$
Then by \eqref{FBSDEs2}, 
\begin{align}\label{widetildex1}
d\widetilde{X}&=(A \widetilde{X}+B^{\trans} \widetilde{u})\, dt+(C \widetilde{X}+D \widetilde{u})^{\trans} dW_t, \\
d\widetilde{Y}&=-[ (Q \widetilde{X}+A \widetilde{Y}+C^{\trans}\widetilde{Z})-n((\be[X]-L)_- - (\be[\widehat{X}]-L)_-)] dt+\widetilde{Z}_t^{\trans} dW_t.\label{widetildey1}
\end{align}
Using It\^o's formula, we have
\begin{align*}
d(\widetilde{X}_t\widetilde{Y}_t)&=B^{\trans} \widetilde{u}\widetilde{Y}\, dt+(C\widetilde{X}+D\widetilde{u})^{\trans} \widetilde{Y} dW_t\nonumber\\
&\quad\;-\widetilde{X}(Q\widetilde{X}-n((\be[X]-L)_- -(\be[\widehat{X}]-L)_-) )dt-\widetilde{Z}^{\trans} dW_t)+\widetilde{Z}_t^{\trans} D\widetilde{u} dt.
\end{align*}
Integrating both sides yields
\begin{align*}
&\quad\;\be[G\widetilde{X}_T^2]+\be\int_0^T Q \widetilde{X}^2_t dt \nonumber\\
&=\int_0^T\be[\widetilde{X}](n(\be[X]-L)_- -n(\be[\widehat{X}]-L)_-)dt-\be\int_0^T\langle R \widetilde{u}_{t}, \widetilde{u}_{t}\rangle dt.\nonumber
\end{align*}
Hence,
\begin{align*}
&\quad\;\be[G\widetilde{X}_T^2]+\be\int_0^T Q \widetilde{X}^2_t dt+\be\int_0^T\langle R \widetilde{u}, \widetilde{u}\rangle dt\le 0,\nonumber
\end{align*}
from which we deduce $\widetilde{u}=0$.
We have $\widetilde{X}=0$ from $\widetilde{X}_0=0$ and \eqref{widetildex1}. This in particular implies $\widetilde{Y}_T=G\widetilde{X}_T=0$. Together with \eqref{widetildey1}, we have $\widetilde{Y}=0$ and $\widetilde{Z}=0$. 

\end{proof}

\begin{Lemma}\label{lemma:optimn}
Suppose that $(X^{n},Y^{n},Z^{n})$ is a solution of the penalized FBSDEs \eqref{FBSDEs2}. 
Then \[u^{n}:=-R^{-1}(BY^{n}+D^{\trans} Z^{n})\] is the optimal control for the unconstrained problem~\eqref{unconstrainedproblem}. And the optimal value is 
\begin{equation}\label{optimalvalue1}
J_n(u^n)={1\over2} {Y}^n_0x+{n\over2}\int_0^T(\be[X_t^n]-L_t)_-L_{t} dt.
\end{equation}
\end{Lemma}
\begin{proof} The proof is similar to that of Theorem \ref{varify1}. We would leave the details to the interested readers. 
\end{proof} 

\subsection{Convergence}
We next show that the solutions of the penalized FBSDEs \eqref{FBSDEs2} have a limit, which turns out to be a solution of the reflected FBSDEs~\eqref{FBSDEs1}. In the following arguments, we may choose a subsequence when necessary. 
Also the constant $M\in \R^+$ might vary from line to line, but does not depend on $n$, $k$ or $t$.

\par
Let $(X^n,Y^n,Z^n, u^{n})$ be given as in Lemma~\ref{lemma:optimn}. Then $u^{n}$ is the optimal control for the unconstrained problem~\eqref{unconstrainedproblem}. 
By Section~\ref{existence}, we conclude that the sequence $u^{n}$ strongly converges to $u^{\infty}$ in the space $L^2_{\cal F}([0,T];\mathbb R)$ and ${X}^n$ converges to ${X}^\infty$ strongly in $C_\mathcal{F}([0,T];\R)$, where $({u}^\infty, {X}^\infty)$ is the optimal pair of the constrained problem~\eqref{cost}. 
Moreover, 
\begin{align} \label{boundedvalue}
0\leq J(u^{n})\le J_n(u^{n})\le J_n(u^{a})=J(u^{a}) \le M,
\end{align} 
where $u^a$ is given in Assumption~\ref{assump1}.

We next show that
$$\mu_t^n:=-n\int_t^T (\be[X_s^n]-L_s)_-ds, \quad t\in [0,T]$$ is a uniformly bounded sequence in $L^\infty ([0, \, T]; \, \R)$ and $Y_{0}^{n}$ is a uniformly bounded sequence in $\R$. 

To this end, let $ \beta\in[x,x+1]$ and $( u^{a}, X^{\beta,a})$ evolve according to equation~\eqref{state} with $ X^{\beta,a}_{0}=\beta$.
By Assumption~\ref{assump1} and monotonicity of SDE, we have $$\be\left[ X_t^{\beta,a}\right]> L_{t},\quad t\in[0,T].$$
Applying Ito's formula to $ X^{\beta, a}_tY_t^n$, we get 
\begin{align*} 
d(X^{\beta, a}Y^{n})
&=Y^{n}(AX^{\beta, a}+B^{\trans} u^{a})\,dt+Y^{n}(CX^{\beta, a}+Du^{a})^{\trans} dW_t \\
&\quad\;-X^{\beta, a}(QX^{n}+AY^n+C^{\trans} Z^n)\,dt+X^{\beta, a}d\mu^{n}_t+X^{\beta, a}(Z^n)^{\trans} dW_t, \\
&\quad\;+(CX^{\beta, a}+Du^{a})^{\trans}Z^n dt\\
&=(-QX^{\beta, a}X^{n}+(u^{a})^{\trans} (BY^{n}+D^{\trans}Z^n))dt+X^{\beta, a}d\mu^{n}_t\\
&\quad\;+(Y^{n}CX^{\beta, a}+Y^{n}Du^{a}+X^{\beta, a}Z^n)^{\trans} dW_t\\
&=(-QX^{\beta, a}X^{n}-(u^{a})^{\trans} Ru^{n}) dt+X^{\beta, a}d\mu^{n}_t\\
&\quad\;+(Y^{n}CX^{\beta, a}+Y^{n}Du^{a}+X^{\beta, a}Z^n)^{\trans} dW_t.
\end{align*}
Integrating on both sides, we have 
\begin{align*} 
\beta Y^n_0&=GX^{\beta,a}_{T} X^n_T+\int_0^T(QX^{\beta, a}X^{n}+(u^{a})^{\trans} Ru^{n})\, dt -\int_0^T X^{\beta,a}d\mu_t^n \\
&\quad\;-\int_0^T(Y^{n}CX^{\beta, a}+Y^{n}Du^{a}+X^{\beta, a}Z^n)^{\trans} dW_t, 
\end{align*} 
and the local martingale is in fact a martingale (see Bismut~\cite[Proposition I-1, p. 387]{Bismut1973}). 
By the elementary inequality $|\left\langle a,b\right\rangle|\leq {1\over 2} |a|^2+{1\over 2}|b|^2$ and \eqref{boundedvalue}, we have 
\begin{align} \label{boundedY}
\abs{\beta Y^n_0+\int_0^T \be[X^{\beta, a}_{t}]d\mu_t^n} &\leq \be[G|X^{\beta, a}_{T} X^n_T|]+\int_0^T\be[Q_{t}|X^{\beta, a}_{t}X^{n}_{t}|+(u^{a}_t)^{\trans} R_{t}u^{n}_t]\, dt\nonumber\\
& \le {1\over 2} \be[G(X^{\beta, a}_{T})^{2}]+ {1\over 2} \int_0^T\be[Q_{t}(X^{\beta, a}_{t})^{2}+(u^{a}_t)^{\trans} R_{t}u^{a}_t]dt\nonumber\\
&\quad\;+ {1\over 2} \be[G(X^{n}_{T})^{2}]+ {1\over 2} \int_{0}^{T}\be[Q_{t}(X^{n}_{t})^{2}+(u^{n}_t)^{\trans} R_{t}u^{n}_t]\, dt\nonumber\\
&\le \frac{1}{2}(M+ J(u^{n})) \leq M, 
\end{align}
where $M$ does not depend on $\beta\in[ x,x+1]$. For the case of $\beta=x$, we have 
\begin{align*} 
\left|xY^n_0+\int_0^T\be[ X^{a}]d\mu_t^n\right| \le M. 
\end{align*}
By \eqref{optimalvalue1} and \eqref{boundedvalue},
\begin{equation*} 
0\leq x{Y}^n_0+ \int_0^T L d\mu^n_t=2J_n(u^{n}) \le M.
\end{equation*}
Comparing the above two inequalities, we get 
\begin{align*} 
\abs{ \int_0^T(\be[ X^{a}]-L)d\mu_t^n }\le M. 
\end{align*}
From Assumption~\ref{assump1} and the monotonicity of $\mu^n$, we see that $\mu^n$ is a uniformly bounded sequence in $L^\infty ([0, \, T]; \, \R)$.
Consequently, choosing $0\neq \beta\in[ x,x+1]$ in \eqref{boundedY}, we see that $Y^n_0$ is a uniformly bounded sequence in $\R$. 

By \eqref{limV}, we have
\begin{align*}
\lim_{n\to\infty}\int_0^T (\be[X^n_t]-L_{t})_- d\mu_{t}^{n}=0.
\end{align*}
Also trivially, 
\begin{align*} 
\int_0^T (\be[X^n_t]-L_{t})_+ d\mu_{t}^{n}=n\int_0^T (\be[X^n_t]-L_{t})_+(\be[X^n_t]-L_{t})_-dt=0,
\end{align*}
so 
\begin{align*}
\lim_{n\to\infty}\int_0^T \abs{\be[X^n_t]-L_{t}} d\mu_{t}^{n}=0.
\end{align*}
As $\mu^n$ is a uniformly bounded sequence, it has a $\star$-weak limit $\mu^{\infty}\in\setmu$. Hence
\begin{align*} 
&\quad\;\int_0^T \abs{\be[X^{\infty}_t]-L_{t}}d\mu^{\infty}_{t} 
=\lim_{n\to\infty}\int_{0}^{T} \abs{\be[X^{\infty}_t]-L_{t}}d\mu^{n}_{t}\\
&\leq\limsup_{n\to\infty}\int_{0}^{T} \abs{\be[X^{\infty}_t]-\be[X^{n}_t]}d\mu^{n}_{t}+\limsup_{n\to\infty}\int_{0}^{T} \abs{\be[X^{n}_t]-L_{t}}d\mu^{n}_{t}\\
&\leq \limsup_{n\to\infty}\max_{t}|\be[X^{\infty}_t]-\be[X^{n}_t]| |\mu^{n}_0|=0,
\end{align*} 
where the last equation is due to the fact that ${X}^n$ converges to ${X}^\infty$ strongly in $C_\mathcal{F}([0,T];\R)$. Therefore, 
\begin{align*} 
&\quad\;\int_0^T (\be[X^{\infty}_t]-L_{t})d\mu^{\infty}_{t}=0. 
\end{align*} 
From the strong convergence of $X^n\to X^\infty$ and equality~\eqref{limV}, we have
\begin{align*}
\int_0^T ((\be[X^{\infty}_t]-L_{t})_-)^2 dt&=\liminf_{n\to\infty}\int_0^T ((\be[X^n_t]-L_{t})_-)^2 dt=0.
\end{align*}
We conclude $\be[X^{\infty}_t]\geq L_{t}$ for all $t\in[0,T]$ by the continuity of $\be[X^{\infty}_{\cdot}]$ and $L_{\cdot}$. 
\par 
Applying the standard estimate for SDE to \eqref{FBSDEs2}, we have 
\begin{gather*}
\be\left[\sup_{t\leq T} |X_t^n|^2+\int_0^T |u^n_t|^2dt\right]\le M,\\
\be\left[\sup_{t\leq T} |X_t^n-X_t^k|^2+\int_0^T |u^n_t-u_t^k|^2 dt\right]\le M, 
\end{gather*}
uniformly for $n$, $k\in\mathbb{N}$.
We will use these estimates for linear SDEs frequently in the subsequent argument without claim. 
\par 
We notice that 
\[\lim_{n,k\to \infty}\be\left[\sup_{t\leq T}|X^n_t-X^k_t|^{2}+\int_{0}^{T}|u^n_t-u^k_t|^{2}dt\right]=0.\]
By H\"{o}lder's inequality, 
\begin{multline*}
\limsup_{n,k\to \infty}\overline{\be}\left[\left(\int_{0}^{T}|u^n_t-u^k_t|^{2}dt\right)^{\frac{p}{2}}\right]\\
\leq \limsup_{n,k\to \infty}\left(\be\left[\int_{0}^{T}|u^n_t-u^k_t|^{2}dt\right]\right)^{\frac{p}{2}}
\left(\be\left[\left(\tfrac{d\overline{\p}}{d \p}\right)^{\frac{2}{2-p}}\right]\right)^{\frac{2-p}{2}}=0, \quad p\in(1,2).
\end{multline*}
Thus,
\begin{align} 
\lim_{n,k\to \infty}\overline{\be}\left[\left(\int_{0}^{T}|u^n_t-u^k_t|^{2}dt\right)^{\frac{p}{2}}\right]=0,\quad p\in(1,2).\label{convergentu}
\end{align}
Similarly, we have
\begin{align}
\lim_{n,k\to \infty}\overline{\be}\left[ \sup_{t\leq T}|X^n_t-X^k_t|^{p} \right]=0,\quad p\in(1,2).
\label{convergentx}
\end{align}

Let $U_{t}=e^{\int_0^t A_rdr}>0$ and $\overline{W}_t=W_t-\int_0^t C_{s}\, ds$. Then $U$ is positive and uniformly bounded. By \eqref{FBSDEs2} 
\begin{equation} 
\begin{cases}
d(UY^n)=-UQX^n\,dt+Ud\mu^n_t+U(Z^n)^{\trans} d\overline{W}_t, \\
Y_T^n=GX_T^n;
\end{cases}
\end{equation}
Integrating yields
\begin{align*}
\int_0^T U(Z^n)^{\trans} d\overline{W}_t=U_TGX_T^n-Y_0^n+\int_0^T UQX^n\,dt-\int_0^T Ud\mu^n_t
\end{align*} 
Because $Y_{0}^{n}$, $\mu^{n}$ and $X^{n}$ are convergent and $U$ is uniformly bounded, by 
\eqref{convergentx}, we have 
\begin{align*}
\lim_{n,k\to \infty}\overline{\be}\left[\left|\int_0^TU_t(Z^n_t-Z^k_t)^{\trans} d\overline{W}_t\right|^{p}\right]=0,\quad p\in(1,2).
\end{align*}
By Doob's martingale inequality,
\begin{align*}
&\quad\;\limsup_{n,k\to \infty}\overline{\be}\left[\sup_{t\leq T}\left|\int_0^t U_s(Z^n_s-Z^k_s)^{\trans} d\overline{W}_s\right|^{p}\right]\\
&\leq \left(\frac{p}{p-1}\right)^{p}\lim_{n,k\to \infty}\overline{\be}\left[\left|\int_0^TU_t(Z^n_t-Z^k_t)^{\trans} d\overline{W}_t\right|^{p}\right]=0,\quad p\in(1,2),
\end{align*}
which gives 
\begin{align*}
&\quad\;\lim_{n,k\to \infty}\overline{\be}\left[\sup_{t\leq T}\left|\int_0^tU_s(Z^n_s-Z^k_s)^{\trans} d\overline{W}_s\right|^{p}\right]=0,\quad p\in(1,2).
\end{align*}
Applying the Burkholder-Davis-Gundy inequality, we conclude 
\begin{align*}
\lim_{n,k\to \infty}\overline{\be}\left[\left(\int_0^TU^2_s|Z_s^n-Z_s^k|^2\,ds\right)^{p\over2}\right]=0,\quad p\in(1,2).
\end{align*}
In turn, by the boundedness of $U$ and H\"{o}lder's inequality
$$
\lim_{n,k\to \infty}\be\left[\left(\int_0^T|Z_s^n-Z_s^k|^2\,ds\right)^{p\over2}\right]=0.
$$ 
As 
$$Y^n=-(B^{\trans}B)^{-1}B^{\trans}(Ru^n+D^{\trans}Z^n),$$
we deduce also
\begin{align} \label{Yconverge}
\lim_{n,k\to \infty}\be\left[\left(\int_0^T|Y_s^n-Y_s^k|^2\,ds\right)^{p\over2}\right]=0.
\end{align}
Now we fix one $p\in(1, 2)$.
Then there exists a unique $(Y^{\infty}, Z^{\infty})\in L^{p}_\cF([0, \, T]; \, \R)\times L^{p}_\cF([0, \, T]; \, \R^m)$ 
such that 
$$
\lim_{n\to \infty}\be\left[\left(\int_0^T|Z_s^n-Z_s^{\infty}|^2\,ds\right)^{p\over2}\right]=0
$$
and 
$$
\lim_{n\to \infty}\be\left[\left(\int_0^T|Y_s^n-Y_s^{\infty}|^2\,ds\right)^{p\over2}\right]=0.
$$

Let $\haty^n=Y^n-\mu^n$. Thanks to \eqref{FBSDEs2} and $\mu_T^n=0$, 
\begin{equation*} 
\begin{cases}
dX^n=(AX^n+B^{\trans}u^n)\,dt+(CX^n+Du^n)^{\trans} dW_t,\\
d\haty^n=-(QX^n+A\haty^n+A\mu^n+C^{\trans} Z^n)\,dt+(Z^n)^{\trans} dW_t, \\
X^n_{0}=x, \quad \haty_T^n=GX_T^n.
\end{cases}
\end{equation*} 
By the standard estimate for BSDE and thanks to the boundedness of the sequence $\mu^n$, 
\begin{align*} 
\be\left[\int_0^T |\haty_t^n|^2\;dt+\int_0^T |Z_t^n|^2\;dt\right]
\le M\Big(\be\left[(GX_T^n)^2\right]+\int_0^T(QX^n+A\mu^n)^2\;dt \Big)\le M. 
\end{align*}
Since $\mu^n$ is a bounded sequence, it yields immediately that
\begin{align*} 
\be\left[\int_0^T |Y_t^n|^2\;dt+\int_0^T |Z_t^n|^2\;dt\right] \le M. 
\end{align*}
Taking lower limits, it follows from Fatou's lemma that
\begin{align*} 
\be\left[\int_0^T |Y_t^{\infty}|^2\;dt+\int_0^T |Z_t^{\infty}|^2\;dt\right]<\infty. 
\end{align*}
So we conclude that 
$(Y^{\infty}, Z^{\infty})\in L^{2}_\cF([0, \, T]; \, \R)\times L^{2}_\cF([0, \, T]; \, \R^m)$. 

Note that $\mu_0^n$ is bounded, we denote its limit (along a subsequence) by $\mu_0$.
Set 
$$\overline{\mu}_t:=Y^\infty_t-Y^\infty_0+\int_0^t (QX_s^\infty+AY_s^\infty+C^{\trans} Z_s^\infty)ds+\mu_0-\int_0^t (Z^\infty_s)^{\trans}dW_s, $$
we have
$$\lim_{n\to \infty}\E{ \left(\int_0^T|\mu_t^n-\overline{\mu}_t|^2\;dt\right)^{p/2}}= 0,\quad p\in (1,2).$$
Let $\varphi\in L_{\cal F}^\infty([0,T])$ and $\Phi=\int_0^t \varphi(s)ds$, then
\begin{align*} 
&\be \int_0^T \varphi(t)\overline{\mu}_t dt=\lim_n \be\left[\int_0^T \varphi(t)\mu^n_t dt\right]=-\lim_n \be\left[\int_0^T \Phi(t)d\mu^n_t \right]
\\
&=-\be\left[\int_0^T \Phi(t)d\mu^\infty_t \right]=\be\left[\int_0^T \varphi(t)\mu^\infty_t dt\right].
\end{align*}
Hence 
$\overline{\mu}=\mu^{\infty},$
and consequently 
$(X^{\infty}, Y^{\infty}, Z^{\infty},\mu^{\infty})$ is a solution to \eqref{FBSDEs1}.

\section{Acknowledgments} 
The three authors would like to thank both referees for their careful reading and helpful comments.

\end{document}